\documentclass[12pt]{article}
\usepackage{geometry}
\usepackage{amsmath}
\usepackage{amsfonts}
\usepackage{amssymb}
\usepackage{amsthm}
\theoremstyle{definition}
\newtheorem{theorem}{Theorem}
\newtheorem{definition}{Definition}
\newtheorem*{corollary}{Corollary}
\newtheorem{lemma}{Lemma}
\begin{document}
\title{Elliptic operators on non-compact manifolds have closed range}
\author{Luther Rinehart}
\maketitle
\begin{abstract}
	We show that a second-order elliptic differential operator $P$, on any manifold $M$, has closed range in $C^\infty(M)$. If $M$ has no compact components, then $P$ is surjective on $C^\infty(M)$. Applications to Helmholtz decomposition are discussed.
\end{abstract}
Most literature on elliptic operators deals with bounded domains and compact manifolds, since in those cases the operators are Fredholm, and the problems of uniqueness as well as existence are well-understood. In contrast this paper deals with existence results in unbounded domains. Specifically, we show that a second-order elliptic operator $P$, on any manifold, compact or not, has closed range in $C^\infty(M)$. Then, for certain manifolds, this combines with the much easier proof that an elliptic operator has dense range, to give the result that the operator is surjective on $C^\infty(M)$.\par
H\"{o}rmander, in section 10.6 of \cite{Hormander}, treats the case of constant coefficient operators on $\mathbb{R}^n$. There they introduced the notion of $P$-convexity, which for constant-coefficient operators is equivalent to $C^\infty$ surjectivity. This proves, for example, that the Laplacian is surjective on $\mathbb{R}^n$, and indeed any open subset of $\mathbb{R}^n$. The method of $P$-convexity can be extended to more general subelliptic operators, as explained by Petersen in section 4.7 of \cite{Petersen}. This approach shows, by functional analytic arguments, that the closed-range property can be reduced to the property of $P$-convexity, an essentially topological property of the domain. Petersen's method is stated for domains in $\mathbb{R}^n$, but is easily adapted to arbitrary manifolds. 
\section{A motivating example}
As a motivating example for existence results in unbounded domains, consider Helmholtz decomposition in elementary vector analysis. It is well known that any vector field in $\mathbb{R}^3$ can be written, non-uniquely, as a sum of a gradient and a curl. For vector fields decaying sufficiently rapidly at infinity, the justification of this fact using, say, the fundamental solution of the Laplacian, is standard textbook material.  However, its proof for general smooth vector fields requires other methods and is less well-known. Consider the following generalization: Let $M$ be an orientable Riemannian manifold, with exterior differential $d$, exterior codifferential $\delta$, and Hodge Laplacian $\Delta=d\delta +\delta d$.
\begin{definition}
	We say $M$ has the \emph{Helmholtz decomposition property} if, given any differential form $\omega$, there exist differential forms $\alpha$ and $\beta$ such that $\omega=d\alpha + \delta\beta$.
\end{definition}
Aside from general usefulness in computations, this property is related to the following facts:
\begin{lemma}\label{decomp}
	If $M$ has the Helmholtz decomposition property, and $\omega=d\alpha + \delta\beta$, then $\alpha$ may be chosen to be coclosed (even coexact), and $\beta$ may be chosen to be closed (even exact).
	\begin{proof}
		Apply decomposition to $\alpha = d\alpha_1 +\delta\alpha_2$ and replace $\alpha$ with $\delta\alpha_2$. Similarly, apply decomposition to $\beta=d\beta_1+\delta\beta_2$, and replace $\beta$ with $d\beta_1$.
	\end{proof}
\end{lemma}
\begin{theorem}\label{Helmholtz}
	$M$ has the Helmholtz decomposition property if and only if $\Delta$ is surjective on differential forms. 
	\begin{proof}
		Clearly, if $\Delta$ is surjective, then given $\omega$, there exists $\phi$ with $\omega=\Delta\phi=d(\delta\phi)+\delta(d\phi)$. \par
		Now suppose $M$ has Helmholtz decomposition. Given $\omega$, we may write 
		$$\omega=d\alpha +\delta\beta$$
		$$\alpha = d\alpha_1 +\delta\alpha_2$$
		$$\beta = d\beta_1+\delta\beta_2$$
		By lemma \ref{decomp}, we can choose $d\alpha_2=0$ and $\delta\beta_1=0$. Then let $\phi=\alpha_2+\beta_1$.
		$$\Delta\phi = (d\delta +\delta d)(\alpha_2+\beta_1)=d\alpha + \delta\beta=\omega$$
	\end{proof}
\end{theorem}
\begin{theorem}
	If $M$ has the Helmholz decomposition property, then the \emph{generalized Maxwell problem} is solvable. That is, given an exact form $\alpha$ and a coexact form $\beta$, there exists a differential form $\omega$ satisfying 
	\begin{equation}
	d\omega = \alpha\qquad\delta\omega=\beta
	\end{equation}
	\begin{proof}
		$\alpha$ is exact, so write $\alpha=d\alpha_1$ where, by Helmholtz decomposition, we may choose $\delta\alpha_1 =0$. Likewise, $\beta$ coexact, so write $\beta=\delta\beta_1$ where, by Helmholtz decomposition, we may choose $d\beta_1=0$. Then  $\omega=\alpha_1+\beta_1$ is the required solution.
	\end{proof}
\end{theorem}

\section{Main results}
Let $M$ be a manifold and $E$ a vector bundle on $M$. Let $P\colon C^\infty(E)\rightarrow C^\infty(E)$ be a second-order elliptic differential operator with smooth coefficients. Give $C^\infty(E)$ its usual Frechet space topology. Let $\mathcal{E}'(E^*) = (C^\infty(E))^*$ denote the space of compactly supported distributional sections. In the following, the bundle $E$ will be suppressed in the notation. A useful fact we will use several times is that if $P$ is second-order elliptic and $Pu=0$, then the support of $u$ must have empty boundary. This first of all gives the easy result:
\begin{theorem}\label{dense}
	Suppose $M$ is connected and not compact. Then $P$ has dense range. 
	\begin{proof}
		 The transpose $P^*$ acting on $\mathcal{E}'$ is also elliptic, so with $M$ being connected and not compact, the equation $P^*u=0$ has no nontrivial compactly supported solutions. Hence $P^*$ is injective, so $P$ has dense range. 
	\end{proof}
\end{theorem}
The rest of this paper will focus on proving that $P$ has closed range. The proof that $P$ has closed range has two halves. The first, functional analytic half, consists of the method laid out in \cite{Petersen} section 4.7. The second half is topological. \par 
Petersen's theorem 7.9 first proves a general subelliptic estimate satisfied by $P$ on compact subsets $K$. This estimate is then used in theorem 7.22 to prove that the image $P(\mathcal{E}'_K)$ of the space of distributions supported in $K$, is weak* closed in $\mathcal{E}'$. The key definition for obtaining a global closed-range theorem is the following.
\begin{definition}
	$M$ is $P$\emph{-convex} if for every compact $K\subseteq M$, there exists a compact $K'\subseteq M$, such that $\forall u\in\mathcal{E}',\ \text{supp}Pu\subseteq K \Rightarrow \text{supp}\,u\subseteq K'$.
\end{definition}
Petersen's corollary 7.24 shows that if $M$ is $P^*$-convex, then $P^*(\mathcal{E}')$ is weak*-closed in $\mathcal{E}'$. By the closed-range theorem for Frechet spaces (\cite{Petersen} p. 33), this is then equivalent to $P(C^\infty)$ being closed in $C^\infty$.\par
If we want to show that $M$ is $P^*$-convex, this amounts to finding, for each compact set $K$, a larger compact set $K'$ such that the components of $M\setminus K'$ are all unbounded, that is, having noncompact closure. This is because the boundary of the support of $u$ cannot intersect the region where $P^*u=0$, so the support of $u$ cannot extend into any unbounded component of $M\setminus K'$. Thus the key topological lemma in this regard is the following one. The proof is thanks to Pierre PC \cite{Pierre} and Aitor Iribar Lopez \cite{Lopez} at mathoverflow.net. It is based on the proof of theorem 23.5 in \cite{Forster}.
\begin{lemma}\label{keylemma}
	Let $X$ be locally compact Hausdorff, connected, and locally connected. Given compact $K\subseteq X$, let $K'$ be the union of $K$ with all the relatively compact components of $X\setminus K$. Then $K'$ is compact.
	\begin{proof}
		Index the connected components of $X\setminus K$ as $\{C_j\}_{j\in J}$, and set $J_0=\{j\in J\mid \overline{C_j}\text{ is compact}\}$. By definition, $K'=K\cup\bigcup_{j\in J_0} C_j$ and $X\setminus K' = \bigcup_{j\notin J_0} C_j$.  As $X$ is locally connected, all the $C_j$ are open, so $K'$ is closed.\par
		As $X$ is locally compact Hausdorff, let $U$ be a relatively compact open neighborhood of $K$. As $X$ is connected, $\forall j,\ C_j\cap U\neq\varnothing$.\par
		Also, $\partial U$ is compact, and $\{C_j\}_{j\in J}$ forms a disjoint open cover of it. Thus there are only finitely many $j_1,...,j_n \in J_0$ such that $C_{j_i}\cap\partial U\neq\varnothing$. For all other $j$, since $C_j$ is connected, we mus have $C_j\subseteq U$. Then,
		$$K' \subseteq U\cup C_{j_1}\cup ...\cup C_{j_n}$$
		Since $K'$ is closed and the right-hand side is relatively compact, it follows that $K'$ is compact.
	\end{proof}
\end{lemma}
\begin{theorem}
	$P$ has closed range.
	\begin{proof}
		First suppose that $M$ is connected. Given a compact $K\subset M$ and constructing $K'$ as in lemma \ref{keylemma}, we see that $M$ is $P^*$-convex. By Petersen's functional analytic argument, it follows that $P$ has closed range. \par
		For non-connected $M$, $C^\infty$ decomposes as the topological direct product of the spaces of smooth sections on the components. Furthermore, differential operators respect this decomposition, so the range of $P$ is the direct product of its ranges on the components, hence closed in the product topology. 
	\end{proof}
\end{theorem}
\begin{corollary}
	If $M$ is connected and not compact, or more generally if $M$ has no compact components, then $P$ is surjective on $C^\infty$.
\end{corollary}
We observe that any failure of surjectivity can be entirely described by the Fredholm theory on the compact components.

\end{document}